\documentstyle[11pt]{article}
\newcommand{\bbn}{\mbox{\boldmath $N$}}
\newcommand{\bbz}{\mbox{\boldmath $Z$}}
\newcommand{\bbc}{\mbox{\boldmath $C$}}
\newcommand{\qed}{\hfill$\Box$}
\newcommand{\Stein}{Steingrimsson}
\newcommand{\maj}{{\sl major}}
\newcommand{\fm}{{\sl flag-major}}
\newcommand{\sgn}{{{\rm Log}_\omega}}
\newcommand{\r}{{r}}
\newcommand{\bg}{{\bar{g}}}
\newcommand{\bp}{{\bar{p}}}
\newcommand{\bq}{{\bar{q}}}
\newcommand{\bx}{{\bar{x}}}
\newcommand{\hg}{{\hat{g}}}
\newcommand{\hG}{{\hat{G}}}
\newcommand{\FZ}{{F_0}}
\newcommand{\FB}{{B}}
\newcommand{\nr}{{n_r}}
\newcommand{\Nri}{{N_r^{(i)}}}
\newcommand{\Nrio}{{N_r^{(i+1)}}}
\newcommand{\Nro}{{N_r^{(1)}}}
\newcommand{\Nrt}{{N_r^{(t)}}}
\newcommand{\pri}{{\pi_r^{(i)}}}
\newcommand{\hpri}{{\hat{\pi}_r^{(i)}}}
\newcommand{\hpro}{{\hat{\pi}_r^{(1)}}}
\newcommand{\hprt}{{\hat{\pi}_r^{(t)}}}
\newenvironment{note}[1]{\par\addvspace{\medskipamount}\noindent
                         {\bf {#1}}\sl
                       }{\par\addvspace{\medskipamount}\rm}

\title{The Flag Major Index and Group Actions on Polynomial Rings}
\author{Ron M. Adin%
\thanks{Department of Mathematics and Computer Science, Bar-Ilan University,
Ramat-Gan 52900, Israel. Email: {\tt radin@math.biu.ac.il}}
$^\ddag$
\and Yuval Roichman%
\thanks{Department of Mathematics and Computer Science, Bar-Ilan University,
Ramat-Gan 52900, Israel. Email: {\tt yuvalr@math.biu.ac.il}}
\thanks{Both authors supported in part by the Israel Science Foundation
and by internal research grants from 
Bar-Ilan University.}}
\date{February 27, 2000}

\begin{document}

\maketitle
\begin {abstract}

A new extension of the major index, defined in terms of Coxeter elements, is introduced.
For the classical Weyl groups of type $B$, it is equidistributed with length.
For more general wreath products
%, as well as for dihedral groups,
it appears in an explicit formula for the Hilbert series of the (diagonal action) 
invariant algebra. 
%An alternative combinatorial interpretation relates this index to the decomposition
%of the coinvariant algebra into irreducible components.

\end{abstract}

\section{Introduction}

\subsection{Outline}

The {\it major index}, $\maj(\pi)$, of a permutation $\pi$ in the symmetric group
$S_n$ is the sum (possibly zero) of all indices $1\le i< n$ for which
$\pi(i)>\pi(i+1)$.
The {\it length} of a permutation $\pi$  is the
minimal number of factors in an expression of $\pi$ as 
a product of the Coxeter generators $(i,i+1), 1\le i< n$.
A fundamental property of the major index is its equidistribution with the
length function [MM]; namely, the number of elements in $S_n$ of a given length $k$ 
is equal to the number of elements having major index $k$.
Bijective proofs and generalizations were given in [F, FS , Ca , GG, Go, Ro3].

\medskip

Candidates for a major index for the classical Weyl groups of type $B$ 
have been suggested by Clarke-Foata [CF1--3], Reiner [Rei1--2],
\Stein\ [Ste], and others.
Unfortunately, unlike the case of the symmetric group, 
the various alternatives are  not equidistributed with
the length function (defined with respect to the Coxeter generators of $B_n$).

In this paper we present a new
definition of the major index for the groups $B_n$, and more generally for wreath
products of the form $C_m\wr S_n$, where $C_m$ is the cyclic group of order $m$.
This major index is shown to be equidistributed with the
length function for the groups of type $B$ (the case $m=2$), and to play a crucial role
in the study of the actions of these groups on polynomial rings.

\medskip

The rest of the paper is organized as follows. The definition of the flag
major index is presented in Subsection 1.2 (and, in more detail,
in Section 2). Main results are surveyed in Subsection 1.3.
Basic properties and equidistribution with length are proved in Section 2.
In Section 3 we give a combinatorial interpretation of the flag major index. 
In Section 4 we study $C_m\wr S_n$-actions on tensor powers of polynomial rings.
The combinatorial interpretation of the flag major index
is then applied to obtain a simple representation of the Hilbert
Series of the diagonal action invariant algebra.

\medskip

\subsection{The Flag Major Index}

The groups $C_m\wr S_n$ are generated
by $n-1$ involutions, $s_1,\ldots, s_{n-1}$,
which satisfy the usual Moore-Coxeter relations of $S_n$,
together with an exceptional generator, $s_0$, of order $m$. 
We consider a different set of generators:
$$
t_i:= \prod_{j=0}^i s_{i-j}\qquad (0 \le i\le n-1).
$$
 These are Coxeter elements in a distinguished flag of parabolic subgroups.
See Section 2 for more details.

An element $\pi\in C_m\wr S_n$ has a unique representation as a product
$$
\pi=t_{n-1}^{k_{n-1}}t_{n-2}^{k_{n-2}}\cdots t_1^{k_1}t_0^{k_0} 
$$
with $0\le k_i< m(i+1)$ $(\forall i)$.
Define the {\it flag major index} of $\pi$ by
$$
\hbox{\fm}(\pi):=\sum\limits_{i=0}^{n-1}k_i.
$$
For $m=1$, this definition gives a new interpretation of a well-known parameter.

\smallskip

\begin{note}
\noindent{\bf Claim 1.} 
For $m=1$ (i.e., for the symmetric group $S_n$) the flag major index coincides with 
the major index.
\end{note}

\noindent See Claim 2.1 below.

\medskip

\subsection{Main Results}

From Claim 1 together with MacMahon's classical result 
it follows that the flag major index  
is equidistributed with length, for $m=1$. This property extends to $m=2$.

\medskip

\begin{note}
\noindent{\bf Theorem 2.} For $m=2$ (i.e., for the hyperoctahedral group $B_n$)
the flag major index is equidistributed with length.
\end{note}

\smallskip

Here ``length" is used in the usual sense, in terms of the Coxeter generators
$s_0,s_1,\ldots,s_{n-1}$. See Theorem 2.2 below.

\smallskip

For $m\ge 3$, the flag major index is no longer equidistributed with length 
(with respect to $s_0,\ldots,s_{n-1}$), but nevertheless 
it does play a central role in the study of naturally defined algebras of polynomials.

\medskip

Let $P_n:=\bbc[x_1,\ldots,x_n]$ be the algebra of polynomials in $n$
indeterminates.
There is a natural action of $G:=C_m\wr S_n$ on $P_n$ (presented explicitly in 
Section 4).  Consider now the tensor power $P_n^{\otimes t}$ with the natural 
{\it tensor action} $\varphi_T$ of $G^t:=G\times\cdots\times G$ ($t$ factors),
and the corresponding {\it diagonal action} of $G$. For more details see Section 4.

The {\it tensor invariant algebra} 
$\hbox{TIA}$ is a subalgebra of the {\it diagonal invariant algebra} 
$\hbox{DIA}$.
Let $F_D(\bq)$, where $\bq=(q_1,...,q_t)$, be the 
multi-variate generating function (Hilbert series) for 
the dimensions of the homogeneous components in $\hbox{DIA}$:
$$
F_D(\bq):= \sum\limits_{n_1,\ldots,n_t\in \bbn}
(\hbox{ dim }_{\bbc}\hbox{ DIA}_{n_1,\ldots,n_t}\hbox{ })
\hbox{ }q_1^{n_1}\cdots q_t^{n_t},
$$
where $\hbox{ DIA}_{n_1,\ldots,n_t}$ is the homogeneous piece of multi-degree
$(n_1,\ldots,n_t)$ in $\hbox{DIA}$. Define similarly
 $F_T(\bq)$ for $\hbox{TIA}$.
Then

\medskip

\begin{note}
\noindent{\bf Theorem 3.} For all $m,n,t\ge 1$,
$$
{F_D(\bq)\over F_T(\bq)} = 
  \sum\limits_{\pi_1\cdots\pi_t=1} \prod_{i=1}^t q_i^{{\fm}(\pi_i)}
$$
where
the sum extends over all $t$-tuples  $(\pi_1,...,\pi_t)$ of elements
in $G=C_m\wr S_n$ such that the product 
$\pi_1\pi_2\cdots\pi_t$ is equal to the identity element.
\end{note}

\medskip

\noindent See Theorem 4.1 below.

\bigskip

\section
{Flag Major Index and Length} 

The groups $C_m\wr S_n$ are generated
by $n-1$ involutions, $s_1,\ldots, s_{n-1}$,
together with an exceptional generator, $s_0$, of order $m$. 
The $n-1$ involutions satisfy the usual Moore-Coxeter relations of $S_n$, 
$$
(s_is_{i+1})^3=1 \qquad (1\le i <n),
$$
$$
(s_is_j)^2=1 \qquad (|i-j|>1),
$$
 while the exceptional generator
satisfies the relations:
$$
(s_0 s_1)^{2m}=1,
$$
$$
s_0 s_i =s_i s_0 \qquad (1< i<n).
$$
For $m=1$, $s_0$ is the identity element.

\medskip

\noindent
Consider now a different set of generators:
$$
t_i:= \prod_{j=0}^i s_{i-j}\qquad (0 \le i\le n-1).
$$
 These are Coxeter elements  [Hu, \S 3.16]
in a distinguished flag of parabolic subgroups
$$
1<G_1<\ldots<G_n=C_m\wr S_n
$$
where $G_i\cong C_m\wr S_i$ is the subgroup of $C_m\wr S_n$ generated by  
$s_0,s_1,\ldots,s_{i-1}$.

An element $\pi\in C_m\wr S_n$ has a unique representation as a product
$$
\pi=t_{n-1}^{k_{n-1}}t_{n-2}^{k_{n-2}}\cdots t_1^{k_1}t_0^{k_0} 
$$
with $0\le k_i< m(i+1)$ $(\forall i)$.
Define the {\it flag major index} of $\pi$ by
$$
\hbox{\fm}(\pi):=\sum\limits_{i=0}^{n-1}k_i.
$$
For $m=1$, this definition gives a new interpretation of a well-known parameter.

\smallskip

\begin{note}
\noindent{\bf Claim 2.1.} 
For $m=1$ (i.e., for the symmetric group $S_n$) the flag major index coincides with 
the major index.
\end{note}

\medskip

\noindent{\bf Proof.} 
Consider the natural action of $S_n$ on the letters $1,\ldots,n$,
where $s_i$ ($1\le i< n$) acts as the transposition $(i,i+1)$.
Here $t_0=s_0$ is the identity permutation, whereas for $1\le r<n$, 
$t_r=s_r\cdots s_1$ acts as the cycle $(r+1,r,\ldots,2,1)$.
The claim may be proved by induction on the exponents $k_i$. 
Obviously, equality holds when $k_1=\ldots=k_{n-1}=0$.
It suffices to prove that for any
 $\pi=t_r^{k_r}\cdots t_1^{k_1}$  
with $1\le r<n$ and $0\le k_r<r$, 
$\maj(t_r\pi)-\maj(\pi)=1$. 
In this case $\pi(j)=j$ for $j>r+1$ and $\pi(r+1)\in\{2,3,\ldots,r+1\}$.
Observe that, for $1\le j\le r+1$, 
$t_r\pi(j)>t_r\pi(j+1)$ if and only if $\pi(j)>\pi(j+1)$,
unless $\pi(j)=1$ or $\pi(j+1)=1$.
If $\pi(j)=1$ (this holds for a unique $1\le j\le r$)
then $t_r\pi(j)=r+1$, so that $\pi(j-1)>\pi(j)<\pi(j+1)$ whereas
$t_r\pi(j-1)<t_r\pi(j)>t_r\pi(j+1)$;
thus the descent at $j-1$ in $\pi$
is replaced by a descent at $j$ in $t_r\pi$.
The special case $j=1$ is similar.

\qed

\medskip

In particular, by MacMahon's classical result, the flag major index  
is equidistributed with length (for $m=1$). This property extends to $m=2$.

\medskip

\begin{note}
\noindent{\bf Theorem 2.2.} 
For $m=2$ (i.e., for the hyperoctahedral group $B_n$),
the flag major index is equidistributed with length.
\end{note}

\smallskip

Here ``length" is used in the usual sense, in terms of the Coxeter generators
$s_0,s_1,\ldots,s_{n-1}$.
Theorem 2.2 can be proved by an explicit bijection.

\medskip

\noindent{\bf Proof.} For $0\le m<2n$
define $r_{n,m}\in B_n$ by
$$
r_{n,m}:= \cases{%
id,                                                  &if $m=0$;\cr
\prod_{j=n-m}^{n-1} s_j,                             &if $0<m\le n$;\cr
\prod_{j=0}^{m-n-1} s_{m-n-j} \prod_{j=0}^{n-1} s_j, &if $n<m<2n$,\cr}
$$
where $id$ is the identity element in $B_n$.

Note that the length of $r_{n,m}$ is $m$.
The set $\{r_{n,m}|0\le m <2n\}$ forms a complete set of 
representatives of minimal length for the left cosets of $B_{n-1}$ in $B_n$.
It follows that every element $\pi\in B_n$ has a unique representation
as a product $\pi=\prod_{i=1}^n r_{n+1-i,m_{n+1-i}}$, 
where $0\le m_j <2j$ for every $j$, and then $\ell(\pi)=\sum_{j=1}^n m_j$.

On the other hand, every element in $B_n$ has a unique representation as
a product of the form
$\prod_{i=1}^n t_{n-i}^{k_{n-i}}$, where $0\le k_j<2(j+1)$ for every $j$.
By definition, 
$\hbox{\fm}(\pi)=\sum_{j=1}^{n} k_{j-1}$.

It follows that the map $\phi:B_n\longrightarrow B_n$ defined by
$$
\phi(\prod_{i=1}^n r_{n+1-i,m_{n+1-i}}):=\prod_{i=1}^n t_{n-i}^{m_{n+1-i}}
$$
is bijective and sends the length function to the flag major index. 

\qed

\medskip

For $m\ge 3$ (when $C_m\wr S_n$ is no longer a Coxeter group), the flag major index 
is no longer equidistributed with length (with respect to $s_0,\ldots,s_{n-1}$); 
nevertheless, it does play a central role in the study of naturally defined algebras 
of polynomials, as will be shown in Section 4.

\medskip

\section{A Combinatorial Interpretation} 
%of the Flag Major Index}

The standard major index for permutations has a natural generalization to any finite
sequence of letters from a linearly ordered alphabet (see, e.g., [F]);
namely, for a finite sequence $a=(a_1,a_2,\ldots,a_n)$ of letters from 
a linearly ordered alphabet, define $\maj(a)$ to be the sum (possibly zero) of 
all indices $1\le i< n$ for which $a_i>a_{i+1}$. 

\smallskip

Let $\omega\in\bbc$ be a primitive $m$-th root of unity.
An element of the wreath product $C_m\wr S_n$ (where $C_m$ is the cyclic group 
of order $m$) may be described as a {\it generalized permutation} 
$\pi = (\pi(1),\pi(2),\ldots,\pi(n))$,
where, for every $i$, $\pi(i)\in \bbc$, 
${\pi(i)\over |\pi(i)|}$ is a power of $\omega$, 
and the sequence of absolute values $|\pi|:=(|\pi(1)|,|\pi(2)|,\ldots,|\pi(n)|)$ 
is a permutation in $S_n$.
In this setting, the involutions $s_i$ ($1\le i< n$) are defined by
$$
s_i(j):=\cases{%
i+1, & if $j=i$;\cr
i,   & if $j=i+1$;\cr
j,   & otherwise,\cr}
$$
whereas the exceptional generator $s_0$ is defined by
$$
s_0(j):=\cases{%
\omega \cdot 1, & if $j=1$;\cr
j,              & otherwise.\cr}
$$
In particular, $C_2\wr S_n$ is the group of {\it signed permutations}, 
also known as the hyperoctahedral group, or the classical Weyl group of type $B$.

\smallskip

Consider now the linearly ordered alphabet
$$
1\cdot \omega^{m-1}<\ldots<n\cdot \omega^{m-1}<\ldots
\ldots<1\cdot \omega^1<\ldots< n\cdot \omega^1<
1\cdot \omega^0<\ldots<n\cdot \omega^0.
$$

\medskip

In this section we prove

\begin{note}
\noindent {\bf Theorem 3.1.} 
 For any $\pi\in C_m\wr S_n$ ,  
$$
\hbox{\fm}(\pi) = m\cdot \maj(\pi) + 
  \sum_{j=0}^{m-1} j\cdot \#\left\{i:{\pi(i)\over |\pi(i)|}=\omega^j\right\}, 
\leqno(3.1)
$$
where $\maj(\pi)$ is defined with respect to the above order.
\end{note}

\bigskip

\noindent To simplify the proof some notations are needed.
For any generalized permutation $\pi\in C_m\wr S_n$ and $1\le i\le n$
define
$$
\sgn\pi(i):=\min\left\{d\ge 0:\omega^d={\pi(i)\over |\pi(i)|}\right\}. 
\leqno(3.2)
$$
Then clearly
$$
\sum_{i=1}^n \sgn\pi(i)=
\sum_{j=0}^{m-1} j\cdot \#\left\{i:{\pi(i)\over |\pi(i)|}= \omega^j\right\}. 
\leqno(3.3)
$$

Denote the right hand side of (3.1) by $\maj_{m,n}(\pi)$. By (3.3),
$$
\maj_{m,n}(\pi)=m\cdot \maj(\pi)+
\sum_{i=1}^n \sgn\pi(i).
\leqno(3.4)
$$

\medskip

\begin{note}
{\bf Lemma 3.2.} For any $\pi\in C_m\wr S_n$ with $\pi(n)\not=\omega^{m-1}$,
$$
\maj_{m,n}(t_{n-1}\pi)-\maj_{m,n}(\pi)=1.
$$
\end{note}

\noindent
{\bf Proof of Lemma 3.2.} 
By definition, 
$$
t_{n-1}(i)=\cases{%
i-1,           & if $i\not=1$;\cr
\omega\cdot n, & if $i=1$.\cr}
$$
Hence, for any $\pi\in C_m\wr S_n$,
$$
t_{n-1}\pi(j)=\cases{%
\omega^{\sgn\pi(j)}\cdot(\pi(j)-1), & if $|\pi(j)|\not=1$;\cr 
\omega^{\sgn\pi(j)+1}\cdot n,       & if $|\pi(j)|=1$.\cr} 
\leqno(3.5)
$$

Let $1\le i_0\le n$ be the unique index for which $|\pi(i_0)|=1$.

\newpage

\noindent{\bf Case (a).} $\sgn\pi(i_0)<m-1$.
\begin{itemize}
\item[]
It follows from (3.5) that if $\sgn\pi(i_0)<m-1$ then 
$$
\maj(t_{n-1}\pi)=\maj(\pi)
$$
with respect to the linear order defined before Theorem 3.1.
Hence, in this case
$$
\maj_{m,n}(t_{n-1}\pi) = 
m\cdot \maj(t_{n-1}\pi) + \sum_{i=1}^n \sgn t_{n-1}\pi(i) =
$$
$$
= m\cdot \maj(\pi) + \sum_{i=1}^n \sgn t_{n-1}\pi(i) =
$$
$$
= m\cdot \maj(\pi) + \sum_{i\not=i_0} \sgn \pi(i)+(\sgn\pi(i_0)+1)
= \maj_{m,n}(\pi)+1.
$$
\end{itemize}

\noindent{\bf Case (b).} $\sgn\pi(i_0)=m-1$.

\begin{itemize}
\item[] 
In this case $\pi(i_0)=1\cdot\omega^{m-1}$ and $t_{n-1}\pi(i_0)=n$. 
By assumption $i_0\not=n$, so that 
$\pi$ has a descent at $i_0-1$ (unless $i_0=1$), 
whereas $t_{n-1}\pi$ has a descent at $i_0$.
Thus 
$$
\maj(t_{n-1}\pi)-\maj(\pi)=1
$$ 
and
$\sgn t_{n-1}\pi(i_0)=0$. Hence, in this case,
$$
\maj_{m,n}(t_{n-1}\pi) 
= m\cdot \maj(t_{n-1}\pi) + \sum_{i=1}^n \sgn t_{n-1}\pi(i) =
$$
$$
= m\cdot (\maj(\pi)+1)+ \sum_{i\not=i_0} \sgn\pi(i) 
= \maj_{m,n}(\pi)+m-(m-1).
$$
\end{itemize}

\qed

\medskip

\noindent
{\bf Proof of Theorem 3.1.} 
By induction on $n$.
Obviously, Theorem 3.1 holds 
in the group $C_m\wr S_1=\langle s_0\rangle$.
Assume that the theorem  holds in the group $C_m\wr S_{n}$, for some $n\ge 1$.
Any element of $C_m\wr S_{n+1}$ has the form
$t_n^{k_n}\pi$, where $\pi\in C_m\wr S_n$ and $0\le k_n< m(n+1)$.
By definition,
$$
\hbox{\fm}(t_n^{k_n}\pi) = k_n+\hbox{\fm}(\pi)
$$
and
$$
\maj_{m,n+1}(\pi)=\maj_{m,n}(\pi).
$$

Hence, by the induction hypothesis, it suffices to prove that
for any $\pi\in C_m\wr S_{n}$ and  $0\le k_n<m(n+1)$
$$
\maj_{m,n+1}(t_{n}^{k_n}\pi)-\maj_{m,n+1}(\pi)=k_n.
$$
This equality readily
follows from iterations of Lemma 3.2, thereby completing the proof.

\qed

\medskip

\section{Diagonal Action on Tensor Powers}

\medskip

Let $P_n:=\bbc[x_1,\ldots,x_n]$ be the algebra of polynomials in $n$
indeterminates.
There is a natural action of $G:=C_m\wr S_n$ on $P_n$, 
$\varphi:G\rightarrow \hbox{Aut}(P_n)$,
defined on generators by 
$$
\varphi(s_0)(x_j)=\cases{%
\omega \cdot x_j, & if $j=1$;\cr
x_j,              & otherwise,\cr} 
\qquad (\omega:=\hbox{exp}(2\pi i/m)\in \bbc)
$$
$$
\varphi(s_i)(x_j)=\cases{%
x_{i+1}, & if $j=i$;\cr
x_i,     & if $j=i+1$;\cr
x_j,     & otherwise,\cr} 
\qquad (1\le i\le n-1)
$$
where each $\varphi(s_i), 0\le i\le n-1$, is extended to an algebra automorphism 
of $P_n$.
Equivalently, in terms of generalized permutations,
$$
\varphi(\pi)(x_j) = {\pi(j)\over |\pi(j)|}\cdot x_{|\pi(j)|} 
\qquad (\forall\pi\in G, 1\le j\le n)
$$
extended multiplicatively to monomials and additively to all of $P_n$.

\medskip

Consider now the tensor power 
$P_n^{\otimes t}:= P_n\otimes \cdots \otimes P_n$ ($t$ factors)
with the natural {\it tensor action} $\varphi_T$ of 
$G^t:=G\times\cdots\times G$ ($t$ factors).
The diagonal embedding
$$
d:G\hookrightarrow G^t
$$
defined by
$$
g\mapsto (g,\ldots,g)\in G^t\qquad (\forall g\in G)
$$
defines the {\it diagonal action} of $G$ on $P_n^{\otimes t}$:
$$
\varphi_D:=\varphi_T\circ d.
$$

\noindent The {\it tensor invariant algebra}
$$
\hbox{TIA}:=
\{\bp\in P_n^{\otimes t}\,|\,\varphi_T(\bg)(\bp)=\bp,\,
\forall \bg\in G^t\}
$$
is a subalgebra of the {\it diagonal invariant algebra} 
$$
\hbox{DIA}:=
\{\bp\in P_n^{\otimes t}\,|\,\varphi_D(g)(\bp)=\bp,\,
\forall g\in G\}.
$$
Note that $\hbox{TIA}=(P_n^G)^{\otimes t}$, where $P_n^G$ is the subalgebra of $P_n$
invariant under $\varphi(G)$.

\medskip

The algebra $P_n^{\otimes t}$ is $\bbn^t$-graded by multi-degree,
where $\bbn:=\{0,1,2,\ldots\}$.
Let $F_D(\bq)$, where $\bq=(q_1,...,q_t)$, be the 
multivariate generating function (Hilbert series) for the dimensions of 
the homogeneous components in $\hbox{DIA}$:
$$
F_D(\bq):= \sum\limits_{n_1,\ldots,n_t\in \bbn}
(\hbox{ dim }_{\bbc}\hbox{ DIA}_{n_1,\ldots,n_t}\hbox{ })
\hbox{ }q_1^{n_1}\cdots q_t^{n_t},
$$
where $\hbox{ DIA}_{n_1,\ldots,n_t}$ is the homogeneous piece of multi-degree
$(n_1,\ldots,n_t)$ in $\hbox{DIA}$. 
Define similarly $F_T(\bq)$ for $\hbox{TIA}$.
Then

\medskip

\begin{note}
\noindent{\bf Theorem 4.1.} 
For all $m,n,t\ge 1$,
$$
{F_D(\bq)\over F_T(\bq)} = 
  \sum\limits_{\pi_1\cdots\pi_t=1} \prod_{i=1}^t q_i^{{\fm}(\pi_i)}
$$
where
the sum extends over all $t$-tuples  $(\pi_1,...,\pi_t)$ of elements
in $G=C_m\wr S_n$ such that the product 
$\pi_1\pi_2\cdots\pi_t$ is equal to the identity element.
\end{note}

\medskip

\section{Proof of Theorem 4.1}

\subsection{Preliminaries}

To prove Theorem 4.1 we need two theorems of Garsia and Gessel.

\medskip

A $t$-{\it partite partition with $n$ parts}
(in the sense of Gordon~[G] and Garsia-Gessel~[GG])
is a sequence $f=(f_1,\ldots,f_t)$ of non-negative-integer valued
functions $f_i:\{1,2,\ldots,n\}\rightarrow \bbn$, $1\le i\le t$,
%$$
%f=\begin{array}{ccccc}
%f_1(1)& f_1(2)&\cdots & f_1(n)\cr
%f_2(1)& f_2(2)&\cdots & f_2(n)\\
%\vdots&\vdots&\vdots&\vdots\\
%f_k(1)& f_k(2)&\cdots & f_k(n)\\
%\end{array}
%$$
satisfying the condition:  
$$
\hbox{If $f_i(j)=f_i(j+1)$ for all $i<i_0$, then $f_{i_0}(j)\ge f_{i_0}(j+1)$
      $\qquad(\forall i_0,j)$.}
$$
In particular, for $i_0=1$: 
$$
f_1(1)\ge f_1(2)\ge\ldots\ge f_1(n)\ge 0,
$$
i.e., $f_1$ is a partition with (at most) $n$ parts.

\smallskip

\noindent{\bf Example.} $n=4, t=2$ ; $f_1=(1,1,0,0)$,
$f_2=(1,0,2,2)$.

\medskip

Let $B_{t,n}$ be the set of all $t$-partite partitions 
$f=(f_1,\ldots,f_t)$ with $n$ parts.
Denote the sum $\sum\limits_{j=1}^n f_i(j)$ by $|f_i|$.
The following theorem was first proved by Garsia and Gessel. 

\smallskip

\begin{note}
\noindent{\bf Theorem GG1} [GG, Theorem 2.2 and Remark 2.2]
$$
\sum_{f\in B_{t,n}} q_1^{|f_1|}\cdots q_t^{|f_t|}=
{\sum_{\pi_1\cdots\pi_t=1} \prod_{i=1}^t q_i^{\maj(\pi_i)}
\over {\prod_{i=1}^t\prod_{j=1}^n (1-q_i^j)}}
$$
where
the sum in the numerator of the right-hand side extends over all $t$-tuples  
$(\pi_1,...,\pi_t)$ of permutations in $ S_n$ such that the product 
$\pi_1\pi_2\cdots\pi_t$ is equal to the identity permutation.
\end{note}

\medskip

Let $n_1,\ldots,n_{\r}$ be non-negative integers such that $\sum_{i=1}^{\r} n_i=n$.
Recall that the {\it $q$-multinomial coefficient} ${n \brack n_1,\ldots,n_{\r}}_q$ 
is defined by:
\[ 
[0]!_q := 1,
\]
\[
[n]!_q := [n-1]!_q \cdot (1 + q + \ldots + q^{n-1}) \qquad (n \ge 1),
\]
\[
{n \brack n_1 \ldots n_{\r}}_q := \frac{[n]!_q}{[n_1]!_q \cdots [n_{\r}]!_q}.
\]
Now let $(N_1,\ldots,N_r)$ be a partition of the set $N:=\{1,\ldots,n\}$. 
For each $1\le i \le \r$ let $\pi_i$ be a permutation on the elements of $N_i$. 
Recall that a permutation $\sigma\in S_n$ is a {\em shuffle} of 
$\pi_1,\pi_2,\ldots,\pi_{\r}$ if, for every $i$, the letters of $N_i$
appear in $\sigma$ in the same order as the corresponding letters appear in $\pi_i$.

\smallskip

\noindent{\bf Example.} 
$N_1=\{1,2,4\}$, $N_2=\{3,5\}$;
$\pi_1 = 241$, $\pi_2 = 35$.  
Here $\sigma=32541$ is a shuffle of $\pi_1$ and $\pi_2$.

\smallskip

\begin{note}
\noindent{\bf Theorem GG2} [GG, Theorem 3.1] 
Let $\Omega(\pi_1,\ldots,\pi_{\r})$ be the collection of 
all shuffles of given permutations
$\pi_1,\pi_2,\ldots,\pi_{\r}$.
Then
$$
\sum_{\sigma\in \Omega(\pi_1,\ldots,\pi_{\r})} q^{\maj(\sigma)}=
{n \brack n_1\ldots n_{\r}}_q q^{\maj(\pi_1)+\ldots +\maj(\pi_{\r})},
$$
where $n_i$ is the number of elements acted upon by $\pi_i$ ($1\le i\le r$).
\end{note}

\bigskip

\subsection{The Case $m=1$}

Before we get to the actual proof of Theorem 4.1,
let us sketch the proof of the case $m=1$ (i.e., $G=S_n$).
This is done as an indication of the structure of the general case;
full details will be given in the actual proof.

The tensor invariant algebra $\hbox{TIA}$ is equal to $(P_n^G)^{\otimes t}$,
and $P_n^G$ is freely generated (as an algebra) by the $n$ elementary symmetric
functions in $n$ indeterminates $x_1,\ldots,x_n$. Thus
$$
F_T(\bq)={1\over \prod_{i=1}^t\prod_{j=1}^n (1-q_i^j)}.
$$
The diagonal invariant algebra $\hbox{DIA}$ is linearly spanned by the
polynomials
$$
\sum\limits_{g\in S_n} \varphi_D(g)(\bx^f)
$$
where $\bx^f$ runs through all the monomials in $P_n^{\otimes t}$. 
Now $\{\bx^f\,|\,f\in B_{t,n}\}$ is a complete set of representatives for the
orbits 
of monomials in $P_n^{\otimes t}$ under the diagonal action $\varphi_D$ of $S_n$.
It follows that a {\it basis} for $\hbox{DIA}$ is
$$
\{\sum\limits_{g\in S_n}\varphi_D(g)(\bx^f)\,|\,f\in B_{t,n}\}.
$$
Therefore, the generating function $F_D(\bq)$ is given by the left-hand side of 
the equation in Theorem GG1.
By Claim 2.1, this proves Theorem 4.1 for $m=1$.
(Theorem~GG2 is not needed in this case.)

\bigskip

In the next subsection we shall construct an explicit basis for 
the diagonal invariant algebra $\hbox{DIA}$, for general $m$. 
In Subsection 5.5 we shall compute its generating function.

\subsection{A Basis for $\hbox{DIA}$}

Consider now the general case of $G=C_m\wr S_n$.
A linear basis for $P_n^{\otimes t}$ consists of the (tensor) monomials
$$
\bx^f:=\bigotimes_{i=1}^t \prod\limits_{j=1}^n x_j^{f_i(j)},
$$ 
where $f_i(j)$ are non-negative integers ($\forall i,j$). Let $F$
denote the set of all such multi-powers $f$. The canonical projection
$\pi:P_n^{\otimes t}\rightarrow \hbox{DIA}$ is defined by
$$
\pi(\bp):= \sum\limits_{g\in G} \varphi_D(g)(\bp) 
\qquad(\forall \bp\in P_n^{\otimes t}),
$$
so that 
$$
\hbox{DIA}= \hbox { span }\{\pi(\bx^f)|f\in F\}.
$$
We are looking for a subset $\FB$ of $F$
such that the set $\{\pi(\bx^f)|f\in \FB\}$ is a basis  for $\hbox{DIA}$.
 
\medskip

First, let
$$
\FZ:=
\{f\in F\,|\,\sum_{i=1}^t f_i(j)\equiv 0\quad(\hbox{mod}\,m)\,\,(\forall j)\}.
$$

\smallskip

\begin{note}
\noindent{\bf Claim 5.1.}
For $f\in F$, 
$$
\pi(\bx^f)\not=0 \Longleftrightarrow f\in \FZ.
$$
\end{note}

\smallskip

\noindent{\bf Proof.}
$$
\varphi_D(s_0)(\bx^f)=\omega^{\alpha(f)}\cdot \bx^f,
$$
where
$$
\alpha(f):=\sum_{i=1}^t f_i(1).
$$
Now, for any $\alpha\in\bbz$:
$$
\sum_{k=0}^{m-1} (\omega^\alpha)^k = \cases{%
m, & if $\alpha \equiv 0\quad(\hbox{mod}\,m)$;\cr
0, & otherwise.\cr}
$$
Therefore, if $C$ is any left coset in $G$ of the subgroup generated by $s_0$,
then $$
\sum_{g\in C}\varphi_D(g)(\bx^f)=0
$$
unless
$$
\sum_{i=1}^t f_i(1)\equiv 0\quad(\hbox{mod}\,m).
$$
Replacing $s_0$ by its $n-1$ conjugates, it follows that 
$$
\pi(\bx^f)\not=0 \Longrightarrow f\in \FZ.
$$

\smallskip

For the converse, let $H$ be the normal commutative subgroup of $G$ generated by 
$s_0$ and its conjugates: $H$ consists of all ``generalized identity permutations''
$h\in G$ satisfying $|h(i)| = i$ $\ (\forall i)$.
Let $\hG$ be the subgroup of $G$ consisting of all ``unsigned permutations''
$\hg\in G$, satisfying $\sgn\hg(i) = 0$ $\ (\forall i)$.
Then $\hG \cong S_n$ and 
%$G = \hG \ltimes H$, 
$G$ is the semidirect product of $\hG$ and $H$, hence 
each $g\in G$ has a unique representation as a product
$$
g=\hg h\qquad(\hg\in\hG, \hbox{ } h\in H).
$$
For any $f\in\FZ$, $\varphi_D(h)(\bx^f)=\bx^f$ for all $h\in H$.
Hence, for $f\in \FZ$, 
$$
\sum\limits_{g\in \hg H} \varphi_D(g)(\bx^f)=|H|\cdot \varphi_D(\hg)(\bx^f)
\qquad(\forall \hg\in\hG).
$$
In other words, the sum over a coset of $H$ in $G$
is a monomial multiplied by 
a positive integer.  It follows that
$$
f\in \FZ \Longrightarrow \pi(\bx^f)\not=0.
$$
\qed

\medskip

For $\bp=\sum\limits_{f\in F} c_f \bx^f\in P_n^{\otimes t}$ (a finite sum), 
define the {\it support}
$$
{\sl supp}(\bp):=\{f\in F\hbox{ }|\hbox{ }c_f\not= 0\}.
$$
The following result is a consequence of the proof of Claim 5.1.

\medskip

\begin{note}
\noindent{\bf Claim 5.2.}
For $f,h\in \FZ$, the following are equivalent:

$(i)$ ${\sl supp}(\pi(\bx^f))\cap {\sl supp}(\pi(\bx^h)) \not= \emptyset$;

$(ii)$ $\pi(\bx^f)=\pi(\bx^h)$;

$(iii)$ $\exists \sigma\in S_n$, such that 
$$
h_i(j)=f_i(\sigma(j))\qquad (\forall i,j).
$$

\end{note}

\medskip

\noindent
Let us now describe explicitly a subset $\FB$ of $\FZ$
such that $\{\pi(\bx^f)|f\in \FB\}$ is a basis for $\hbox{DIA}$.
This subset will, of course, be a complete set of representatives for 
the orbits of all monomials $\{\bx^f|f\in\FZ\}$ under the action of $S_n$ 
described in Claim~5.2(iii).

Intuitively, we classify the exponents $f_i(j)$ according to their vector of 
residues modulo $m$, and attach a $t$-partite partition to each possible vector
of residues.
Formally, define
$$
R:=\{(r_1,\ldots,r_t)\in \bbz^t\,|\,0\le r_i< m\,\,(\forall i)\hbox{\ and\ }
\sum_{i=1}^t r_i\equiv 0\,\,(\hbox{mod}\,m)\}
$$
and choose an arbitrary linear order $\le_R$ on $R$.

\smallskip

Define a bijection
$$
\theta:\FZ\longrightarrow F\times R^n
$$
by
$$
\theta(f):=(h,r),
$$
where $h=(h_i(j))_{i,j}\in F$ and $r=(r_i(j))_{i,j}\in R^n$ are defined as the
quotients and remainders, respectively, obtained when the entries of
$f=(f_i(j))_{i,j}\in \FZ$ are divided by $m$: 
$$
f_i(j)=m\cdot h_i(j)+r_i(j) \qquad (\forall i,j).
$$ 

\noindent
Now let $\FB$ be the set of all $f\in \FZ$ such that $\theta(f)=(h,r)$
satisfies:
\begin{itemize}
\item[(i)] 
$r_*(1)\ge_R \ldots \ge_R r_*(n)$, where 
$r_*(j):=(r_1(j),\ldots,r_t(j))\in R$ ($\forall j$);
\item[(ii)] 
If $r_*(j)=r_*(j+1)$, and also $h_i(j)=h_i(j+1)$ for all $i<i_0$,
then 
$$
h_{i_0}(j)\ge h_{i_0}(j+1).
$$
\end{itemize}

\smallskip

We can interpret $f\in \FB$ as follows: 
For each $r_*\in R$ there is a range (possibly empty) of indices $j_1< j\le j_2$
for which $r_*(j)=r_*$, and
a sequence of vectors $(h_*(j))_{j=j_1+1}^{j_2}$
which forms a $t$-partite partition with $j_2 - j_1$ parts
(as defined in Subsection~5.1 above).
The total number of vectors (for all $r_*\in R$) is $n$.

\medskip

Clearly, $\FB$ is a complete system of representatives for the orbits of $\FZ$ 
under the action of $S_n$ defined in Claim~5.2(iii). 
By Claims~5.1 and~5.2 we conclude

\medskip

\begin{note}
\noindent{\bf Lemma 5.3.}
The set 
$$
\{\pi(\bx^f)\,|\,f\in \FB\}
$$
is a homogeneous basis for $\hbox{DIA}$.
\end{note}

\medskip

\begin{note}
\noindent{\bf Corollary 5.4.}
The generating function for $\hbox{DIA}$ is
$$
F_D(\bq) = \sum_{f\in\FB} q_1^{|f_1|} \cdots q_t^{|f_t|} =
\sum_{(\nr)} \left[\prod_{r\in R} (q_1^{r_1} \cdots q_t^{r_t})^{\nr}
  F_{t,\nr}(q_1^m,\ldots,q_t^m)\right], 
$$
where $F_{t,n}(q_1,\ldots,q_t)$ is the generating function for $t$-partite partitions
described in Theorem GG1, and the sum is over all partitions of $n$ into non-negative 
integers $(\nr)$ indexed by $r\in R$.
\end{note}

\medskip 

\noindent{\bf Example.} $m=t=2$.
\begin{itemize}
\item[]  
In this case $R=\{(1,1),(0,0)\}$, so that
$f=(f_1(j),f_2(j))_{1\le j\le n}$ belongs to $\FB$ if and only if
there exists an integer $0\le k\le n$ such that: 
\begin{itemize}
\item[(i)] 
For $1\le j\le k$, $f_1(j)$ and $f_2(j)$ are both odd (i.e., $r_*(j)=(1,1)$); 
and for $k<j\le n$ they are both even (i.e., $r_*(j)=(0,0)$).
\item[(ii)] 
The vector pairs $(h_1(j),h_2(j))_{1\le j\le k}$ and
$(h_1(j),h_2(j))_{k< j\le n}$ are both $2$-partite partitions.
\end{itemize}
\end{itemize}

\medskip

\medskip

\subsection{Constructing Sequences of Generalized Permutations}

In this subsection we construct a bijection
between $t$-tuples of generalized permutations in $G$ with product equal to the
identity element, and data consisting of numbers, sets and permutations 
defined below. This bijection will allow us to
apply Theorem GG2 in the computation of the generating function,
to be carried out in Subsection 5.5.

\medskip

Every sequence $(\pi_1,\ldots,\pi_t)\in G^t$ such that $\pi_t \cdots \pi_1 = 1$
may be constructed using the following steps:

\begin{description}
\item[{\bf Step 1.}] 
Choose non-negative integers $(\nr)_{r\in R}$ such that 
$$\sum_{r\in R} \nr = n.$$
\item[{\bf Step 2.}] 
For each $1\le i\le t$ choose a set-partition $(\Nri)_{r\in R}$ of
$N:=\{1,\ldots,n\}$ such that 
$$\#\Nri = \nr\qquad(\forall r\in R),$$
where $\# S$ denotes the size of the set $S$.
\item[{\bf Step 3.}] 
For each $1\le i\le t-1$ and $r\in R$ choose a bijective function (``permutation'')
$\hpri:\Nri\to\Nrio$ (essentially, $\hpri\in S_{\nr}$);
and define $\hprt:\Nrt\to\Nro$ such that $\hprt \cdots \hpro = 1$. 
\item[{\bf Step 4.}] 
For each $1\le i\le t$ define a permutation $\hat{\pi}_i \in S_n$ and a generalized
permutation $\pi_i \in G$ as follows:
if $1\le j\le n$ and  $j \in \Nri$ then
$$\hat{\pi}_i(x_j) := \hpri(x_j)$$
and
$$\pi_i(x_j) := \omega^{r_i}\cdot \hpri(x_j).$$
\end{description}

Conversely, let $(\pi_1,\ldots,\pi_t)\in G^t$ satisfy $\pi_t \cdots \pi_1 = 1$.
For each $1\le j\le n$ and $1\le i\le t$ let
$$\r_i(j) := \sgn\pi_i(|\pi_{i-1} \cdots \pi_1(j)|).$$
Then 
$$r_*(j) := (r_1(j),\ldots,r_t(j)) \in R,$$
since $\pi_t \cdots \pi_1 = 1$.

Now define $(\nr)_{r\in R}$ and $(\Nri)_{r\in R}$ by:
$$\Nro := \{1\le j\le n \,|\, r_*(j) = r\}\qquad(\forall r\in R),$$
$$\nr := \# \Nro\qquad(\forall r\in R),$$
$$\Nri := \{|\pi_{i-1} \cdots \pi_1(j)| \,:\, j\in\Nro\}
  \qquad(\forall r\in R, 2\le i\le t),$$
and $$\pri := \pi_i |_{\Nri}
  \qquad(\forall r\in R, 1\le i\le t).$$
To sum up, there is a bijection
$$(\pi_1,\ldots,\pi_t) \leftrightarrow \left((\nr),(\Nri),(\hpri)\right)$$
between $t$-tuples of generalized permutations in $G$ with product equal to the
identity element, and data consisting of numbers, sets and permutations as above.

\smallskip

\subsection{Computation of the Generating Function}

In this subsection we prove the claim of Theorem~4.1, namely that
$$
\frac{F_D(\bq)}{F_T(\bq)} = 
\sum_{{\pi_1,\ldots,\pi_t \in G} \atop {\pi_t \cdots \pi_1 = 1}}
  \prod_{i=1}^t q_i^{\fm(\pi_i)}. \leqno(5.1)
$$

\smallskip

\noindent
From Theorem~3.1 and from the bijection in the previous subsection it follows that 
the right-hand side of (5.1) is equal to
$$
\hbox{RHS}=\sum_{{\pi_1,\ldots,\pi_t\in G} \atop {\pi_t \cdots \pi_1 = 1}}
  \prod_{i=1}^t q_i^{\fm(\pi_i)} = \leqno (5.2)
$$
$$
=\sum_{(\nr)} \sum_{(\Nri)} \sum_{(\hpri)}
  \prod_{i=1}^t q_i^{m\cdot \maj(\pi_i) + \sum_{r\in R} \nr r_i} =
$$
$$
= 
\sum_{(\nr)} \left[\prod_{r\in R} (q_1^{r_1} \cdots q_t^{r_t})^{\nr}
  \sum_{(\Nri)} \sum_{(\hpri)} \prod_{i=1}^t q_i^{m\cdot \maj(\pi_i)}\right].
$$

\noindent Now, $\hbox{TIA}=(P_n^G)^{\otimes t}$, and $P_n^G$ consists of all
the symmetric polynomials in $x_1^m,\dots,x_n^m$. Thus
$$
F_T(\bq) = \frac{1}{\prod_{i=1}^t \prod_{j=1}^n (1-q_i^{mj})}, \leqno (5.3)
$$
Combining (5.3) with Corollary~5.4 (from the end of Subsection 5.3),
we obtain that the left-hand side of (5.1) is equal to
$$
\hbox{LHS}=\frac{F_D(\bq)}{F_T(\bq)}=
\prod_{i=1}^t \prod_{j=1}^n (1-q_i^{mj})\cdot
\sum_{(\nr)} \left[\prod_{r\in R} (q_1^{r_1} \cdots q_t^{r_t})^{\nr}
  F_{t,\nr}(q_1^m,\ldots,q_t^m)\right].
\leqno(5.4)
$$
%$$
%=\sum_{(\nr)} \prod_{r\in R}(q_1^{r_1} \cdots q_t^{r_t})^{\nr}
%\cdot F_{t,\nr}(q_1^m,\ldots,q_t^m)\cdot
%\prod_{i=1}^t \prod_{j=1}^n (1-q_i^{mj}).
%$$
Comparing (5.2) with (5.4) we conclude that
it suffices to show that, for every choice of $(\nr)_{r\in R}$,
$$ 
\prod_{i=1}^t \prod_{j=1}^n (1-q_i^{mj}) \cdot 
  \prod_{r\in R} F_{t,\nr}(q_1^m, \ldots, q_t^m) =
\sum_{(\Nri)} \sum_{(\hpri)} \prod_{i=1}^t q_i^{m \cdot \maj(\pi_i)}.
\leqno(5.5)
$$

Theorem~GG1 above gives an explicit expression for $F_{t,\nr}$:
$$
F_{t,\nr}(q_1^m, \ldots, q_t^m) = 
\left[\prod_{i=1}^t \prod_{j=1}^{\nr} (1-q_i^{mj})\right]^{-1} \cdot 
%  \sum_{{\hpro, \ldots, \hprt \in S_{\nr}} \atop {\hprt \cdots \hpro = 1}}
  \sum_{\hprt \cdots \hpro = 1}
  \prod_{i=1}^t q_i^{m\cdot \maj(\hpri)}, 
$$
where $\hpri\in S_{n_r}$ are unsigned permutations.

Denoting $q:=q_i^m$, the definition of $q$-multinomial coefficients gives
$$
\prod_{j=1}^n (1-q^j) \cdot 
  \left[\prod_{r\in R} \prod_{j=1}^{\nr} (1-q^j)\right]^{-1} =
{n \brack (\nr)_{r\in R}}_q,
$$
so that the left-hand side of~(5.5) is equal to
$$
\prod_{i=1}^t {n \brack (\nr)_{r\in R}}_{q_i^m} \cdot 
  \sum_{(\hpri)} \prod_{r\in R} \prod_{i=1}^t q_i^{m\cdot \maj(\hpri)}.
$$
Here the sum is over all choices of $\hpri\in S_{\nr}$ ($r\in R$, $1\le i\le t$)
such that $\hprt \cdots \hpro = 1$ ($\forall r$). 

Thus, all we need to prove is that, for every choice of $(\nr)$ and $(\hpri)$,
$$
\prod_{i=1}^t  {n \brack (\nr)}_{q_i^m} \cdot 
  \prod_{r\in R} \prod_{i=1}^t q_i^{m\cdot \maj(\hpri)} =
\sum_{(\Nri)} \prod_{i=1}^t q_i^{m\cdot \maj(\pi_i)}. \leqno(5.6)
$$
Since the choice of the partition $(\Nri)_{r\in R}$ of $N$ can be made
independently for each value of $i$, we can ``interchange'' the sum and product in
the right-hand side of~(5.6), and thus it suffices to show that (again, denoting
$q:=q_i^m$):
$$
{n \brack (\nr)}_q \cdot \prod_{r\in R} q^{\maj(\hpri)} =
\sum_{(\Nri)_{r\in R}} q^{\maj(\pi_i)}\qquad(\forall i, (\nr), (\hpri)).
\leqno (5.7)
$$
Now this amounts exactly to the statement of Theorem~GG2, since the different 
choices of $(\Nri)_{r\in R}$ correspond to the various shuffles of the permutations
$(\hpri)_{r\in R}$, each yielding a different $\pi_i$.

\qed

%\newpage

If $t=2$, a somewhat simpler argument may be given.  We demonstrate it in the special
case $m=t=2$.

\medskip

\noindent{\bf Example:} $m=t=2$ .

\medskip

In this case, by Corollary~5.4, Theorem GG1 and (5.4) (using simplified notation) 
$$
{F_D(\bar q)\over F_T(\bar q)}
=\prod\limits_{i=1}^2\prod\limits_{j=1}^n (1-q_i^{2j})\cdot \left\{
\sum\limits_{k=0}^n (q_1 q_2)^k \cdot 
{\sum\limits_{\pi_1\in S_{k}} q_1^{2\cdot\maj(\pi_1)}
                              q_2^{2\cdot\maj(\pi_1^{-1})}
\over \prod\limits_{i=1}^2\prod\limits_{j=1}^{k} (1-q_i^{2j})}
\right\}
\cdot
$$
$$
\cdot
\left\{
{\sum\limits_{\pi_1\in S_{n-k}} q_1^{2\cdot \maj(\pi_1)}
                                q_2^{2\cdot \maj(\pi_1^{-1})}
\over \prod\limits_{i=1}^2\prod\limits_{j=1}^{n-k} (1-q_i^{2j})}
\right\}=
$$
$$
=
\sum\limits_{k=0}^n q_1^k q_2^k \cdot 
\prod\limits_{i=1}^2
{n\brack k}_{q_i^2} \Sigma_{k,i}
%\left(\sum\limits_{\pi_1\in S_{k}} q_1^{2\cdot \maj(\pi_1)}
%                                   q_2^{2\cdot \maj(\pi_1^{-1})}\right)
%\cdot
%$$
%$$
%\cdot
%\left(\sum\limits_{\pi_1\in S_{n-k}} q_1^{2\cdot \maj(\pi_1)}
%                                   q_2^{2\cdot \maj(\pi_1^{-1})}\right)
%=
$$
where
$$
\Sigma_{k,i} = 
\left(\sum\limits_{\pi_1\in S_{k}} q_1^{2\cdot \maj(\pi_1)}
                                   q_2^{2\cdot \maj(\pi_1^{-1})}\right)
\cdot
\left(\sum\limits_{\pi_1\in S_{n-k}} q_1^{2\cdot \maj(\pi_1)}
                                     q_2^{2\cdot \maj(\pi_1^{-1})}\right)
=
$$
$$
%\sum\limits_{k=0}^n q_1^k q_2^k \cdot 
%\prod\limits_{i=1}^2
%{n\brack k}_{q_i^2}
\sum\limits_{\pi_1\in S_{k}\wedge \pi_2\in S_{n-k}} 
q_1^{2\cdot(\maj(\pi_1)+\maj(\pi_2))}
q_2^{2\cdot(\maj(\pi_1^{-1})+\maj(\pi_2^{-1}))}.
$$
Now, observe that every signed permutation
in $B_n=C_2\wr S_n$ may be constructed in the following way:
\begin{description}
\item[{\bf Step 1.}] 
Choose a number $0\le k\le n$.
\item[{\bf Step 2.}] 
Choose $k$ digits from the set $\{1,\ldots,n\}$ and mark them ``negative" 
(the rest will be ``positive").
\item[{\bf Step 3.}] 
Choose two permutations: $\pi_1\in S_k$ (on the ``negative" digits), 
and $\pi_2\in S_{n-k}$ (on the ``positive" ones).
\item[{\bf Step 4.}] 
Choose a shuffle of $\pi_1$ and $\pi_2$.
\end{description}
It should be noted that the major index (in the sense defined at 
the beginning of Section 3) of the resulting signed permutation
is independent of the choice in Step 2.
In fact, it is exactly the major index of the shuffle (chosen in Step 4) 
of $\pi_1$ and $\pi_2$.
%\noindent
On the other hand, the major index of the inverse of the resulting signed 
permutation does not depend on the choice in Step 4. 
In fact, it is exactly the major index of the shuffle (defined by Step 2) 
of $\pi_1^{-1}$ and $\pi_2^{-1}$.

Combining these facts with Theorem GG2, it follows that 
for any choosen integer $0\le k\le n$ and any given 
pair of permutations $\pi_1\in S_k$ and $\pi_2\in S_{n-k}$
$$
\prod\limits_{i=1}^2
{n\brack k}_{q_i^2} q_1^{2\cdot(\maj(\pi_1)+\maj(\pi_2))}
q_2^{2\cdot(\maj(\pi_1^{-1})+\maj(\pi_2^{-1}))}= $$ $$ =\sum\limits_{\sigma\in
B(\pi_1,\pi_2)} q_1^{2\cdot \maj(\sigma)} q_2^{2\cdot \maj(\sigma^{-1})}, 
$$
where $B(\pi_1,\pi_2)$ is the set of all signed permutations
constructed by choosing $\pi_1$ and $\pi_2$ at the third step.

Thus,
$$ 
\prod\limits_{i=1}^2
{n\brack k}_{q_i^2} \sum\limits_{\pi_1\in S_{k}\wedge \pi_2\in S_{n-k}}
q_1^{2\cdot(\maj(\pi_1)+\maj(\pi_2))}
q_2^{2\cdot(\maj(\pi_1^{-1})+\maj(\pi_2^{-1}))}= $$ $$ =\sum\limits_{\sigma\in
B_n(k)} q_1^{2\cdot \maj(\sigma)} q_2^{2\cdot \maj(\sigma^{-1})}, 
$$ 
where
$B_n(k):=\{\sigma\in B_n| \sigma \hbox{ has $k$ ``negative" digits}\}$.  
We conclude that
$$
{F_D(\bq)\over F_T(\bq)}= \sum\limits_{k=0}^n q_1^k q_2^k \cdot
\sum\limits_{\sigma\in B_n(k)} q_1^{2\cdot \maj(\sigma)} q_2^{2\cdot
\maj(\sigma^{-1})}= $$ $$ = \sum\limits_{\sigma\in B_n} q_1^{2\cdot
\maj(\sigma)+k(\sigma)} q_2^{2\cdot \maj(\sigma^{-1})+k(\sigma)}, 
$$ 
where $k(\sigma)$ is the number of ``negative" digits in $\sigma$. Note that
$k(\sigma)=k(\sigma^{-1})$. Theorem 3.1 completes the proof of the desired
result (for $m=t=2$).

\bigskip

\section{\bf Final Remarks} 
\begin{itemize}
\item[\bf 1.] 
The flag major index may be defined on dihedral groups in an analogous way 
(with respect to the Coxeter generators).
An exact analogue of Theorem 4.1 may be derived. This will be proved elsewhere.
\item[\bf 2.]
Using Theorem 3.1 it is possible to connect the flag major index with multiplicities
of irreducible representations in homogeneous components of the 
coinvariant algebra of the group $C_m\wr S_n$. These multiplicities were calculated
by Stembridge [Stem] and involve major indices of skew standard 
Young tableaux. 
For more details see [AR, Section 5].
\end{itemize}

\bigskip

\noindent
{\bf Acknowledgments.} Thanks to Richard Stanley, Dominique Foata, 
Christian Krattenthaler, Igor Pak and Victor Reiner for their comments.

\end{document}